\input amstex
\input epsf
\newdimen\unit
\unit=0.025in
\def\plot#1 #2 #3 {\rlap{\kern#1\unit\raise#2\unit\hbox{$#3$}}}
\magnification=\magstep1
\baselineskip=13pt
\documentstyle{amsppt}
\vsize=8.7truein
\NoRunningHeads
\def\EE{\operatorname{\bold E}}
\def\PP{\operatorname{\bold P}}
\def\haf{\operatorname{haf}}
\topmatter
\title Random Weighting, Asymptotic Counting,
 and Inverse Isoperimetry \endtitle
\author Alexander Barvinok and  Alex Samorodnitsky \endauthor
\address Department of Mathematics, University of Michigan, Ann Arbor,
MI 48109-1109, USA \endaddress
\email barvinok$\@$umich.edu  \endemail
\address Department of Computer Science, Hebrew University of 
Jerusalem, Givat Ram Campus, Jerusalem 91904, Israel \endaddress
\email salex$\@$cs.huji.ac.il  \endemail
\date June 2003 \enddate
\thanks The research of the first
author was partially supported by NSF Grant DMS 9734138.
The research of the second author was partially supported by ISF Grant 
039-7165.
\endthanks
\abstract 
For a family $X$ of $k$-subsets of the set $\{1, \ldots, n\}$, let 
$|X|$ be the cardinality of $X$ and let 
$\Gamma(X, \mu)$ be the expected maximum weight of a subset from
$X$ when the weights of $1, \ldots, n$ are chosen independently at 
random from a symmetric probability distribution $\mu$ on ${\Bbb R}$. 
We consider the inverse isoperimetric problem of finding $\mu$ for which 
$\Gamma(X, \mu)$ gives the best estimate of $\ln |X|$. 
We prove that the optimal choice of $\mu$ is the logistic distribution,
in which case $\Gamma(X, \mu)$ provides an asymptotically tight estimate 
of $\ln |X|$ as $k^{-1} \ln |X|$ grows.
Since in many important cases $\Gamma(X, \mu)$ can be easily computed,
we obtain computationally efficient approximation algorithms for a 
variety of counting problems. 
Given $\mu$, we describe families $X$ of a given cardinality 
with the minimum value of $\Gamma(X, \mu)$, thus extending and 
sharpening various isoperimetric inequalities in the Boolean cube.
\endabstract
\endtopmatter
\document

\head 1. Introduction \endhead

Let $X$ be a family of $k$-subsets of the set $\{1, \ldots, n\}$.
Geometrically, we think of $X$ as a set of points 
$x=(\xi_1, \ldots, \xi_n)$ in the {\it Hamming sphere} of radius 
$k$
$$\xi_1 + \ldots + \xi_n=k \quad \text{where} \quad \xi_i \in \{0,1\}
\quad \text{for} \quad i=1, \ldots, n.$$ 
We also consider general families $X$ of subsets of
$\{1, \ldots, n\}$, which we view as sets $X \subset \{0,1\}^n$ 
of points in the Boolean cube. 

Let us fix a Borel probability measure $\mu$ in ${\Bbb R}$. 
We require $\mu$ to be symmetric, that is, $\mu(A)=\mu(-A)$ for any 
Borel set $A \subset {\Bbb R}$, and to have finite variance.

In this paper, we relate two quantities associated with $X$. The first 
quantity is the cardinality $|X|$ of $X$. 
The second quantity
$\Gamma(X, \mu)$ is defined as follows. 
Let us fix a measure $\mu$ as above and let $\gamma_1, \ldots, \gamma_n$ 
be independent random variables having the distribution $\mu$.  
Then 
$$\Gamma(X, \mu)=\EE \max_{x \in X} \sum_{i \in x} \gamma_i.$$
In words: we sample weights of $1, \ldots, n$ independently at random 
from the distribution $\mu$, define the weight of a subset 
$x \in X$ as the sum of the weights of its elements and let 
$\Gamma(X, \mu)$ be the expected maximum weight of a subset from $X$.

Often, when the choice of $\mu$ is clear from the context or not important, 
we write simply $\Gamma(X)$.

It is easy to see that $\Gamma(X)$ is well defined, 
that $\Gamma(X)=0$ if $X$ consists of a single point 
(recall that $\mu$ is symmetric) and that 
$\Gamma(X) \geq \Gamma(Y)$ provided $Y \subset X$. 
Thus, in a sense, $\Gamma(X)$ measures how large $X$ is.
In some respects, $\Gamma(X)$ behaves rather like $\ln |X|$.
For example, if $X \subset \{0, 1\}^n$ and 
$Y \subset \{0, 1\}^m$, we can define the direct product 
$X\times Y \subset \{0, 1\}^{m+n}$. In this case, 
$|X \times Y|=|X| \cdot |Y|$ and $\Gamma(X \times Y)=\Gamma(X)+\Gamma(Y)$. 

Our goal can be stated (somewhat vaguely) as follows:
\bigskip
\noindent
{\bf (1.1) Problem.} Find a measure $\mu$ for which $\Gamma(X, \mu)$ gives
the best estimate of $\ln |X|$. 
\bigskip
Our motivation comes from problems of efficient combinatorial counting. 
For many interesting 
families $X$, given a set $\gamma_1, \ldots, \gamma_n$ of 
weights, we can easily find the maximum weight of a subset $x \in X$
using well-known optimization algorithms. The value of $\Gamma(X, \mu)$ 
can be efficiently computed through 
averaging of several sample maxima for randomly chosen weights 
$\gamma_1, \ldots, \gamma_n$. At the same time, counting elements in 
$X$ can be a hard and interesting problem. 
Thus, for such families, $\Gamma(X, \mu)$ provides a quick estimate 
for $\ln |X|$.
We give some examples in Section 2, where we also argue that the problems of
optimization (computing $\Gamma(X, \mu)$) and counting (computing 
$\ln |X|$) are asymptotically equivalent.
\subhead (1.2) The logistic measure \endsubhead
Let $X$ be a non-empty family of $k$-subsets of $\{1, \ldots, n\}$.
One of our main results is that there are measures $\mu$ for 
which $\Gamma(X, \mu)$ gives an asymptotically 
tight estimate for $\ln |X|$ provided $\ln |X|$ grows faster than 
a linear function of $k$. We obtain the best estimates
when $\mu=\mu_0$ is the {\it logistic measure} with density
$${1 \over e^{\gamma} +e^{-\gamma} +2} \quad \text{for} \quad 
\gamma \in {\Bbb R}.$$
We prove that for any $\alpha>1$ there exists $\beta=\beta(\alpha)>0$ 
such that 
$$\beta \Gamma(X) \leq \ln |X| \leq \Gamma(X) \quad 
\text{provided} \quad |X| \geq \alpha^k$$ 
and 
$$\beta(\alpha) \longrightarrow 1 \quad \text{as} \quad 
\alpha \longrightarrow +\infty.$$
Moreover, we prove that for $t=k^{-1} \Gamma(X)$ we have
$$t-\ln t -1 \leq k^{-1} \ln |X| \leq t$$
for all sufficiently large $t$.
Note that the bounds do not depend on $n$ at all.

The existence of such measures $\mu$ seems to contradict 
some basic geometric intuition. If we fix the cardinality $|X|$ of 
a set $X$ in the Hamming sphere, we would expect $\Gamma(X)$ to be 
large if $X$ is ``random'' and small if $X$ is tightly packed.
It turns out, however, that there are measures that manage to ignore,
to some extent, the difference between dense and sparse sets. 
In Sections 4 and 5 we prove some general asymptotically tight bounds
which allow one to obtain similar estimates 
for a variety of measures $\mu$. For example, the measure with 
density $|\gamma| e^{-|\gamma|}/2$ also guarantees the 
asymptotic equivalence of $\ln |X|$ and $\Gamma(X)$.

We prove that the logistic distribution is, in a well-defined sense, 
optimal among 
all distributions $\mu$ for which $\ln |X| \leq \Gamma(X, \mu)$ 
for all non-empty $X \subset \{0, 1\}^n$:
given a lower bound for $\Gamma(X, \mu)$, we get the best lower
bound for $\ln |X|$ when $\mu$ is the logistic distribution. 
 
In addition, we prove that the
logistic distribution has an interesting extremal property:
the inequality $\ln |X| \leq \Gamma(X)$ turns into equality if $X$ is 
a face (subcube) of the Boolean cube $\{0, 1\}^n$. 

We state our results in Section 3.
\bigskip
The problems we are dealing with have obvious connections to some 
central questions in probability and combinatorics, such as discrete 
isoperimetric inequalities (cf. \cite{ABS98}, \cite{Le91}, and \cite{T95}) 
and estimates of the supremum of a stochastic process, see \cite{T94}.
In particular, in \cite{T94}, M. Talagrand considers 
the functional $\Gamma(X, \mu)$, where $X$ is a family of 
subsets of the set $\{1, \ldots, n\}$ and $\mu$ is the symmetric 
exponential distribution with density $e^{-|\gamma|}/2$. He proves that 
$\ln |X| \leq c \Gamma(X)$ for some absolute constant $c$,
see also \cite{La97}. In Section 7, we prove that the optimal value of 
the constant is $c=2\ln 2$ (the equality is obtained when $X$ is a 
face of the Boolean cube $\{0,1\}^n$). We also prove that 
$\ln |X| \leq \Gamma(X)+ k \ln 2$ provided $X$ lies in the 
Hamming ball of radius $k$ (the inequality is asymptotically sharp).

\subhead (1.3) Isoperimetric inequalities \endsubhead
Suppose that $\mu$ is the Bernoulli measure:
$$\mu\{1\}=\mu\{-1\}={1 \over 2}.$$
This case was studied in our paper \cite{BS01}.
It turns out that $\Gamma(X)$ has a simple geometric interpretation:
the value of $0.5 n- \Gamma(X)$ is the 
average Hamming distance from a point $x$ in the Boolean cube 
$\{0, 1\}^n$ to the subset $X \subset \{0,1\}^n$.
The classical isoperimetric inequality in the Boolean cube, 
Harper's Theorem (see \cite{Le91}), implies
that among all sets $X$ of a given cardinality, the smallest value of 
$\Gamma(X)$ is attained when $X$ is the sphere in the Hamming metric.
More precisely, let us fix $0 < \alpha < \ln 2$. Then there exists 
$\beta=\beta(\alpha)$, $0< \beta <1/2$, such that if $Y_n$ is the Hamming
sphere of radius $\beta n +o(n)$ in $\{0,1\}^n$ then 
we have $\ln |Y_n|=\alpha n +o(n)$ and for any set $X_n \subset \{0,1\}^n$
with $\ln |X_n|=\alpha n +o(n)$, 
we have $\Gamma(Y_n) \leq \Gamma(X_n) +o(n)$. 
We determine $\beta$ from the equation 
$$\beta \ln {1 \over \beta} + (1-\beta) \ln {1 \over 1-\beta}=\alpha$$ and 
note that $\Gamma(Y_n)=\beta n + o(n)$.

In Section 8, we construct sets $Y_n$ with asymptotically the 
smallest value of $\Gamma(Y_n)$ for an {\it arbitrary} 
symmetric probability measure $\mu$ with finite variance. 
It is no longer true that $Y_n$ 
is a Hamming sphere in $\{0,1\}^n$. For example, 
if $\mu\{1\}=\mu\{-1\}=\mu\{0\}=1/3$ then $Y_n$ has to be 
the direct product of two Hamming spheres. It turns out that 
for any symmetric $\mu$ with finite variance $Y_n$ 
can be chosen to be the direct product of at most two Hamming spheres.   
More precisely, let us fix a symmetric probability measure $\mu$ and 
a number $0< \alpha < \ln 2$. Then we construct numbers 
$\lambda_i=\lambda_i(\alpha, \mu) \geq 0$ and 
$0 \leq \beta_i=\beta_i(\alpha, \mu) \leq 1/2$ for $i=1,2$, such that  
$\lambda_1 + \lambda_2=1$ and the following holds: if $Y_n$ is 
the direct product of the Hamming sphere of radius 
$ \beta_1 n +o(n)$ in the Boolean cube of dimensions
$\lambda_1 n +o(n)$ and the Hamming sphere of 
radius $\beta_2 n+o(n) $ in the Boolean cube of 
dimension $\lambda_2 n +o(n)$ (so that $Y_n$ is a subset of 
the Boolean cube of dimension $n$) then $\ln |Y_n|=\alpha n +o(n)$ 
and $\Gamma(Y_n) \leq \Gamma(X_n) +o(n)$ for any 
set $X_n \subset \{0,1\}^n$ such that $\ln |X_n|=\alpha n +o(n)$.

\subhead (1.4) The inverse isoperimetric problem \endsubhead
It turns out that for
the inequality $\ln |X| \leq c \Gamma(X, \mu)$ to hold
with some constant $c=c(\mu)$ for any non-empty family $X$ 
of $k$-subsets of $\{1, \ldots, n\}$, the measure 
$\mu$ has to have an at least exponential tail, that is,
for the cumulative distribution function $F$ of $\mu$ we must have 
$1-F(t) \geq e^{-at}$ for some $a>0$ and all sufficiently large $t$.  
On the other hand, for the lower bound of $k^{-1}\ln |X|$ in terms of
$k^{-1}\Gamma(X, \mu)$ to be non-trivial 
(other than 0), $\mu$ has to have an at most exponential tail.
Thus for the solution of Problem 1.1, which we call the {\it inverse}
isoperimetric problem, we are interested in measures with exponential 
tails.

\head 2. Applications to Combinatorial Counting \endhead

This research is a continuation of \cite{B97} and \cite{BS01},
where the idea to use optimization algorithms for counting problems was 
developed. 

First, we discuss how to compute $\Gamma(X)$
for many interesting families of subsets. 

Let us assume that the family $X$ of subsets of 
$\{1, \ldots, n\}$ is given by its 
{\it Optimization Oracle}.
\specialhead (2.1) Optimization Oracle \endspecialhead
\bigskip
\noindent{\bf Input:} Real vector $c=(\gamma_1, \ldots, \gamma_n)$ 
\medskip
\noindent{\bf Output:} Real number 
$$w(X,c)=\max_{x \in X} \sum_{i \in x} \gamma_i.$$
\bigskip
Thus, we input real weights of the elements $1, \ldots, n$ and 
output the maximum weight $w(X,c)$ of a subset $x \in X$ in this 
weighting. As is discussed in \cite{B97} and 
\cite{BS01}, for many interesting families
$X$ Optimization Oracle 2.1 can be easily constructed. We provide 
two examples below. 
\example{(2.2) Bases in matroids} Let $A$ be a $k \times n$
 matrix of rank $k$ over a field ${\Bbb F}$. We assume that $k<n$. 
Let $X=X(A)$ be the set of all 
$k$-subsets $x$ of $\{1, \ldots, n\}$ such that the columns of $A$ 
indexed by the elements of $x$ are linearly independent. 
Thus $X$ is the set of all non-zero $k \times k$ minors of $A$, or,
in other words, the set of {\it bases} of the matroid represented 
by $A$.
It is an interesting and apparently hard problem to compute or 
to approximate the cardinality of $X$, cf. \cite{JS97}.

On the other hand, it is 
very easy to construct the Optimization Oracle for $X$. Indeed,
given real weights $\gamma_1, \ldots, \gamma_n$, we construct 
a linearly independent set $a_{i_1}, \ldots, a_{i_k}$ of columns 
of the largest total weight one-by-one.
First, we 
choose $a_{i_1}$ to be a non-zero column of $A$ with the 
largest possible weight $\gamma_{i_1}$. Then
we choose $a_{i_2}$ to be a column of the maximum possible 
weight such that $a_{i_1}$ and $a_{i_2}$ are linearly independent.
We proceed as above, and finally select $a_{i_k}$ to be a 
column of the maximum possible weight such that 
$a_{i_1}, \ldots, a_{i_k}$ are linearly independent; cf., for 
example, Chapter 12 of \cite{PS98} for ``greedy algorithms''.
Particular cases of this problem include counting forests and spanning 
subgraphs in a given graph. 

Let $A$ and $B$ be $k \times n$ matrices of rank $k<n$ 
and let $X$ be the set of all $k$-subsets $x$ of $\{1, \ldots, n\}$ such that the columns of 
$A$ indexed by the elements of $x$ are linearly independent and the 
columns of $B$ indexed by the elements of $x$ are linearly independent.
Then there exists a much more complicated than above,
but still polynomial time algorithm, which, given weights $\gamma_1,
\ldots, \gamma_n$, computes the largest weight of a subset $x$ from
$X$, see Chapter 12 of \cite{PS98}. 
\endexample 
\example{(2.3) Perfect matchings in graphs}
Let $G$ be a graph with $2k$ vertices and $n$ edges. A 
collection of $k$ pairwise disjoint edges in $G$ is called 
a {\it perfect matching} (known to physicists as a {\it dimer cover}). 
It is a hard and 
interesting problem to count perfect
matchings in a given graph, see \cite{JS97}. Recently, using the 
Markov chain approach,
M. Jerrum, A. Sinclair and E. Vigoda constructed a polynomial time 
approximation algorithm to count perfect matchings in 
a given bipartite graph \cite{JSV01}, but for general 
graphs no such algorithms are known.

There is a classical $O(n^3)$ algorithm for finding 
a perfect matching of the maximum weight in any given edge-weighted
graph, see Section 11.3 of \cite{PS98}, so Oracle 2.1 is 
readily available.  
\endexample

For any set $X$ given by 
its Optimization Oracle 2.1, the value of $\Gamma(X)$ can be 
well approximated by the sample mean of a moderate size.

\specialhead (2.4) Algorithm for computing $\Gamma(X, \mu)$ 
\endspecialhead 
\bigskip
\noindent{\bf Input:} A family $X$ of subsets 
of $\{1, \ldots, n\}$ given by its Optimization 
Oracle 2.1;
\medskip
\noindent{\bf Output:} A number $w$ approximating 
$\Gamma(X, \mu)$; 
\medskip
\noindent{\bf Algorithm:} Choose a positive integer $m$ 
(see Section 2.5 for details).
Sample independently $m$ random vectors $c_i$ from the 
product measure
$\mu^{\otimes n}$ in ${\Bbb R}^n$. For each vector 
$c_i$, using Optimization Oracle 2.1, compute the maximum 
weight $w(X, c_i)$ of a subset from $X$. Output 
$$w={1 \over m} \sum_{i=1}^m w(X, c_i).$$
\subhead (2.5) Choosing the number of samples $m$ \endsubhead 
Let us consider the output 
$$w=w(X; c_1, \ldots, c_m)$$ of Algorithm 2.4 as a random variable 
on the space 
$${\Bbb R}^{nm}=\underbrace{{\Bbb R}^n \oplus \ldots \oplus
{\Bbb R}^n}_{m \text{\ times}}$$  
endowed with the product measure $\mu^{\otimes mn}$.
Clearly, the expectation of $w$ is $\Gamma(X, \mu)$.

Let $D=\EE(\gamma^2)$ be the variance of $\mu$. Using the 
estimates
$$\Bigl( \sum_{i \in x} \gamma_i \Bigr)^2 
\leq \Bigl( \sum_{i=1}^n |\gamma_i| \Bigr)^2 \leq 
n \sum_{i=1}^n \gamma_i^2 \quad \text{for} \quad x \subset \{1, \ldots, n\},
$$
we conclude that the variance of $w$ does not exceed 
$n^2D/m$. Therefore, by Chebyshev's inequality, 
for the output $w$ to satisfy $|w-\Gamma(X, \mu)| \leq \epsilon$ 
with probability at least $2/3$, we can choose 
$m=\lceil 3 \epsilon^{-2} n^2 D \rceil$.

As usual, to achieve a higher probability $1-\delta$ of success,
we can run the algorithm $O(\ln \delta^{-1})$ times and then 
find the median of the computed estimates.

For many measures $\mu$ the bound for $m$ can be essentially improved.
In particular, we are interested in the case of the logistic measure 
$\mu$ with density $(2+e^{\gamma}+ e^{-\gamma})^{-1}$.
To obtain the desired estimate we use a concentration property of the 
symmetric exponential measure $\nu$ with density $e^{-|\gamma|}/2$,
see Section 4.5 of \cite{Led01}.

Let us define  
$$\psi(\gamma)=\cases 
\gamma-\ln\bigl(2-e^{\gamma}\bigr) &\text{if} \  \gamma \leq 0 \\
\gamma+\ln\bigl(2-e^{-\gamma}\bigr) &\text{if} \ \gamma >0 \endcases$$
and
$$\Psi(c)=\bigl(\psi(\gamma_1), \ldots, \psi(\gamma_n)\bigr) \quad 
\text{for} \quad c=(\gamma_1, \ldots, \gamma_n).$$
Then $\psi(\gamma)$ has the logistic distribution $\mu$ if $\gamma$ has 
the exponential distribution $\nu$. Thus we can write 
$$\Gamma(X, \mu)=\EE w\bigl(X; \Psi(c_1), \ldots, \Psi(c_m)\bigr),$$
where vectors $(c_1, \ldots, c_m)$ are sampled from the 
exponential distribution $\nu^{\otimes mn}$ in ${\Bbb R}^{nm}$. 
If $X$ is a family of 
$k$-subsets then the Lipschitz coefficient of 
$$f(c_1, \ldots, c_m)=w\bigl(X; \Psi(c_1), \ldots, \Psi(c_m) \bigr)$$
with respect to the $\ell^2$ metric of ${\Bbb R}^{nm}$ does not exceed 
$2 \sqrt{k/m}$ while the Lipschitz coefficient with respect to the 
$\ell^1$ metric does not exceed $2/m$. 
Applying Proposition 4.18 of \cite{Led01},
we conclude that for the output $w$ of Algorithm 2.4 to satisfy 
$|w-\Gamma(X, \mu)| \leq \epsilon$ with probability at least $2/3$,
we can choose $m=O(k\epsilon^{-2})$. Note that the choice of $m$ is 
independent of the size $n$ of the ground set.

We observe that it is easy to sample a random 
weight $\gamma$ from the logistic distribution provided sampling 
from the uniform distribution on the interval 
$[0,1]$ is available (which is the case for many computer 
packages). Indeed, if $\xi$ is uniformly distributed on the interval
$[0,1]$, then $\gamma=\ln \xi - \ln (1-\xi)$ has the logistic 
distribution. 

Our numerical experiments suggest that the choice of 
$m=O(1)$ (for example, $m=5$) is good enough and that in many cases 
$m=1$ suffices.
\subhead (2.6) Counting with multiplicities \endsubhead
Suppose that every element $i$ of the ground set $\{1, \ldots, n\}$ 
has a positive integer {\it multiplicity} $q_i$. Let $X$ be a 
family of $k$-subsets of $\{1, \ldots, n\}$ and let 
$$p_X(q_1, \ldots, q_n)=\sum_{x \in X} \prod_{i \in x} q_i.$$
It may be of interest to compute or approximate $p_X$. 

For instance, let $A=(a_{ij})$ be a $2k \times 2k$ symmetric 
matrix of non-negative integers $a_{ij}$. Let us construct 
an (undirected) graph $G$ on $2k$ vertices $\{1, \ldots, 2k\}$ 
where the vertices $i$ and $j$ are connected by an edge if and 
only if $a_{ij}>0$. We identify the edges of $G$ with 
the set $\{1, \ldots, n\}$. 
Let $X$ be the set of all perfect matchings in $G$ identified 
with a family of $k$-subsets of $\{1, \ldots, n\}$, see Example 2.3.
If we assign multiplicities $a_{ij}$ to the edges of $G$, then 
the value of $p_X(a_{ij})$ is called the {\it hafnian} $\haf A$ of 
$A$, a polynomial of a considerable interest which generalizes permanent.

Computing $p_X(q_1, \ldots, q_n)$ is reduced to counting in the 
following straightforward way. Let $N=q_1 + \ldots + q_n$ and 
let us view the set $\{1, \ldots, N\}$ as the multiset
consisting of $q_1$ copies of $1$, $q_2$ copies of $2$, $\ldots$, $q_n$ copies of $n$.
Let us construct a family $Y$ of $k$-subsets of $\{1, \ldots, N\}$ 
as follows: for each $k$-subset $x \in X$ we construct 
$\prod_{i \in x} q_i$ $k$-subsets $y \in Y$ by replacing every 
$i \in x$ by any of its $q_i$ copies. It is clear that 
$|Y|=p_X(q_1, \ldots, q_n)$. 

To construct Optimization Oracle 2.1 for $Y$, we apply the 
oracle for $X$ with the input $c=(\gamma_1, \ldots, \gamma_n)$,
where $\gamma_i$ is the maximum of $q_i$ weights assigned to 
the $q_i$ copies of $i$. Moreover, Algorithm 2.4 is easily modified
for computing $\Gamma(Y, \mu)$ instead of $\Gamma(X, \mu)$.
We still work with the underlying family $X$, but instead of sampling 
weights from the distribution $\mu$, we sample the $i$-th weight
$\gamma_i$ from the distribution $\mu_{q_i}$ of the maximum of $q_i$ 
independent 
random variables with the distribution $\mu$
(note that $\mu_{q_i}$ is not symmetric for $q_i>1$). Thus, if $\mu$ is 
the logistic distribution, to sample $\gamma_i$, we sample $\xi$ from 
the uniform distribution on $[0,1]$ 
and let $\gamma_i=-\ln\bigl(\xi^{-1/q_i} -1 \bigr)$.
Luckily, for the logistic distribution the required number $m$ of calls
to Oracle 2.1 does not depend 
on the size of the ground set, hence we use the same number 
$m$ of calls whether we consider counting with or 
without multiplicities.     
\bigskip
In \cite{BS01} we discuss how our 
approach fits within the general framework of the Monte Carlo method.
The estimates we get are not nearly as precise as 
those obtained by the Markov chain based Monte 
Carlo Method (see, for example, \cite{JS97}), but supply a non-trivial
information and are easily computed for a wide variety of problems. 
Even for the much-studied problem of 
counting perfect matchings in general (non-bipartite) graphs our approach 
produces new theoretical results. 
For some of the problems, such as counting bases 
in the intersection of two general matroids (see Example 2.2), 
our estimates seem to be the only ones that can be 
efficiently computed at the moment. If $X$ is a family of $k$-subsets 
of $\{1, \ldots, n\}$ and $|X|=e^{k \lambda}$ for some 
$\lambda=\lambda(X)$ then, in polynomial time, we estimate 
$\lambda(X)$ within a constant multiplicative factor as long as $\lambda(X)$ is separated from 0 and
all sufficiently large $\lambda(X)$ are estimated with  
an additive error of $1+\ln \lambda(X)$, see Section 3.
Similar estimates hold for counting with multiplicities of Section 2.6.
On the other hand, the Markov chain approach, if successful, 
allows one to estimate the cardinality $|X|$ within any prescribed relative error. We note that for truly large problems the correct scale is logarithmic
because $|X|$ can be prohibitively large to deal with. The 
Markov chain approach relies on the local structure of $X$ (needed 
for ``rapid mixing''), whereas our method uses some global structure
(the ability to optimize on $X$ efficiently).
\subhead (2.7) 
Asymptotic equivalence of counting and optimization \endsubhead 
One can view the optimization functional 
$\max_{x \in X} \sum_{i \in x} \gamma_i$ as the ``tropical version''
of the polynomial $p_X(q_1, \ldots, q_n)$ of Section 2.6: we get the 
former if we replace ``$+$'' with ``$\max$'' and product with sum 
in the latter. Thus our results establish a weak asymptotic 
equivalence of the counting and optimization problems: if we 
can optimize, we can estimate $\ln p_X$ with a relative error 
which approaches 0 as $k^{-1} \ln p_X$ grows. Vice versa, if we 
can approximate $\ln p_X$, we can optimize (at least approximately):
choosing 
$q_i(t)=2^{t \gamma_i}$, we get
$$\lim_{t \longrightarrow +\infty} t^{-1} \log_2 
p_X\bigl(q_1(t), \ldots, q_n(t)\bigr)=
\max_{x \in X} \sum_{i \in x} \gamma_i.$$  
A. Yong \cite{Y03} implemented our algorithms for some counting problems,
such as estimating the number of forests in a given graph, 
computing the permanent and the
hafnian of a given non-negative integer matrix and 
performed a number of numerical experiments. The algorithm produces 
the upper and lower bounds for the logarithm of the cardinality 
of the family in question, see Section 3. The upper bound is 
attained when the family is ``tightly packed'' as a subset of the 
Boolean cube whereas the lower bound is attained on sparse families.
It appears that 
there is some metric structure inherent to various families of 
combinatorially defined sets. For example, when we applied our methods 
to estimate the logarithm of 
the number of spanning trees in a given connected graph,
the exact value (which can be easily computed by the matrix-tree 
formula) turns out to be very close to the upper bound obtained by our
algorithm. Informally, spanning trees appear to be ``tightly 
packed''. On the other hand, when we 
estimated the logarithm of the 
number of perfect matchings in a graph, the true value 
(when we were able to find it by other methods) 
seems to lie close to the middle point between the upper and lower 
bounds. 
   
\head 3. The Logistic Measure: Results \endhead

Let us choose $\mu_0$ with density
$${1 \over e^{\gamma} + e^{-\gamma} +2} \quad 
\text{for} \quad \gamma \in {\Bbb R}.$$
The cumulative distribution function $F$ of $\mu$ is given 
by
$$F(\gamma)={1 \over 1+e^{-\gamma}} \quad \text{for} \quad 
\gamma \in {\Bbb R}.$$
The variance of $\mu_0$ is $\pi^2/3$ \cite{M85}.
Our first main result is as follows.
\proclaim{(3.1) Theorem} 
\roster
\item For every non-empty set $X \subset \{0, 1\}^n$, we have 
$$\ln |X| \leq \Gamma(X);$$
\item
Let 
$$h(t)=\sup_{0 \leq \delta < 1} \Bigl(\delta t +
\ln {\sin \pi \delta \over \pi \delta} \Bigr) \quad 
\text{for} \quad t \geq 0.$$
Then $h(t)$ is a convex increasing function and for 
any non-empty family $X$ of $k$-subsets of $\{1, \ldots, n\}$, we
have 
$$h(t) \leq  k^{-1} \ln |X| \quad \text{where} \quad t=k^{-1} \Gamma(X).$$
\endroster
\endproclaim 
From the expansion
$$\ln {\sin \pi x \over \pi x} =-{\pi^2 \over 6} x^2 +O(x^4) \quad \text{for} 
\quad x \approx 0,$$
we deduce that 
$$h(t)={3 \over 2 \pi^2} t^2 +O(t^4) \quad 
\text{for} \quad t \approx 0$$
(we substitute $\delta=(3t/\pi^2)$).
From the expansion
$$\ln { \sin \pi (1-x) \over \pi(1-x)} =\ln x + x +O(x^2)
\quad \text{for} \quad x \approx 0,$$
we deduce that
$$h(t) \geq  t-\ln t -1 \quad \text{as} \quad t \longrightarrow +\infty$$
(we substitute $\delta=1-t^{-1}$).

A Maple plot of $h(t)$ is shown on Figure 1 below.
$$\hbox to 2.0 in{
\plot -30 0 {\epsffile{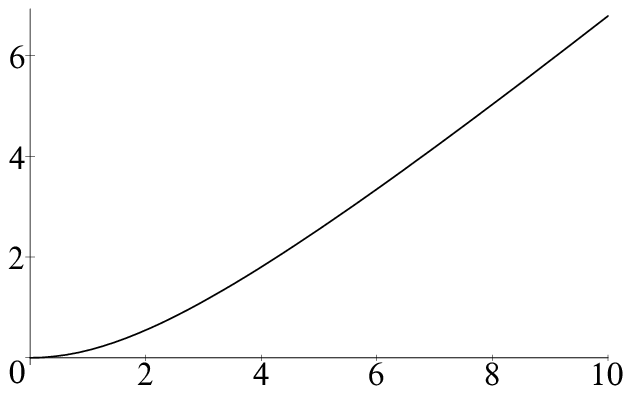}} 
\plot 70 80 {h(t)} 
\plot 85 20 {t}
\plot 35 10 {\bold{Figure \ 1}} 
\hfil}$$
We obtain the following corollary.
\proclaim{(3.2) Corollary} For any $\alpha>1$ there exists 
$\beta=\beta(\alpha)>0$ such that for any non-empty family 
$X$ of $k$-subsets of $\{1, \ldots, n\}$ with $|X| \geq \alpha^k$
we have 
$$\beta \Gamma(X) \leq \ln |X| \leq \Gamma(X).$$
Moreover, 
$$\beta(\alpha) \longrightarrow 1 \quad \text{as} \quad 
\alpha \longrightarrow +\infty.$$
\endproclaim 
\demo{Proof} From Part (1) of Theorem 3.1, we have 
$k^{-1} \Gamma(X) \geq \ln \alpha$. Since $h(t)$ is convex,
we have $h(t) \geq \beta t$ for some $\beta=\beta(\alpha)>0$ 
and all $t \geq \ln \alpha$. The asymptotics of $\beta(\alpha)$ 
as $\alpha \longrightarrow +\infty$ follows from the 
asymptotics of $h(t)$ as $t \longrightarrow +\infty$. 
{\hfill \hfill \hfill} \qed
\enddemo
Thus, using the logistic distribution allows us to estimate 
$\ln |X|$ within a constant factor and the approximation factor 
approaches 1 as $k^{-1} \ln |X|$ grows.

We note that the bound $\ln |X| \leq \Gamma(X)$ is sharp.
For example, if $X$ is an $m$-dimensional face of the Boolean cube 
then $\ln |X|=m \ln 2$ and one can show that $\Gamma(X)=m \ln 2$ 
as well. Indeed, because $\Gamma(X)$ is invariant under 
coordinate permutations, we may assume that $X$ consists 
of the points $(\xi_1, \ldots, \xi_m, 0, \ldots, 0)$, where 
$\xi_i \in \{0, 1\}$ for $i=1, \ldots, m$. The set $X$ can be 
written as the Minkowski sum $X=X_1+ \ldots +X_m$, where 
$X_i$ consists of the origin and the $i$-th basis vector $e_i$.
Hence $\Gamma(X)=m\Gamma(X_1)$ (cf. Section 4.1) and 
$\Gamma(X_1)$ is computed directly as 
$$\Gamma(X_1)=\int_0^{+\infty} {x \over e^x +e^{-x} +2} \ dx=
\ln 2$$
(we substitute $e^x=y$ and then integrate by parts).

It turns out that the logistic measure is optimal in a well-defined 
sense.
\proclaim{(3.3) Theorem} Let ${\Cal M}$ be the set of all measures $\mu$ 
such that 
$$\ln |X| \leq \Gamma(X, \mu)$$
for any non-empty family $X$ of $k$-subsets of $\{1, \ldots, n\}$, any 
$n \geq 1$, and any $1 \leq k \leq n$.

For a measure $\mu \in {\Cal M}$ and a number $t>0$, let 
$c(t, \mu)$ be the infimum of $k^{-1} \ln |X|$ taken over all $n \geq 1$,
all $1 \leq k \leq n$, and all
non-empty families $X$ of $k$-subsets $\{1, \ldots, n\}$ such that 
$k^{-1} \Gamma(X, \mu) \geq t$.
Then for all $t>0$
$$c(t, \mu) \leq c(t, \mu_0),$$
where $\mu_0$ is the logistic distribution.
\endproclaim
\subhead (3.4) Discussion \endsubhead Unless $\mu$ is concentrated 
in $0$, for $X=\{0,1\}$ we have $\Gamma(X, \mu)=c\ln 2$ for some $c>0$ and hence 
$\Gamma(X, \mu)=c \ln |X|$ if $X$ is a face of the Boolean cube 
$\{0,1\}^n$, cf. Section 4.1. As we are looking for the best measure $\mu$
in Problem 1.1,
it is only natural to assume that $\Gamma(X, \mu) \geq c_1 \ln |X|$ for 
all $X \subset \{0,1\}^n$, which, after scaling, becomes 
$\Gamma(X, \mu) \geq \ln |X|$. This explains the definition 
of ${\Cal M}$.

Let us choose $\mu \in {\Cal M}$. 
Then any upper bound for $\Gamma(X, \mu)$ 
is automatically an upper bound for $\ln |X|$. The function $c(t, \mu)$
measures the quality of the lower bound estimate for $\ln |X|$ given a lower
bound for $\Gamma(X, \mu)$. 

Incidentally, it follows from our proof 
that the logistic measure is the measure of the smallest 
variance in ${\Cal M}$.
\bigskip
We prove Theorems 3.1 and 3.3 in Section 6.

\head 4. General Estimates: the Upper Bound \endhead

It is convenient to think about families $X$  
geometrically, as subsets of the Boolean cube 
$\{0, 1\}^n \subset {\Bbb R}^n$.  Let us fix a symmetric probability 
measure $\mu$ on ${\Bbb R}$ with finite variance and let 
$\mu^{\otimes n}$ be the product measure on ${\Bbb R}^n$. 
For a finite set $X \subset {\Bbb R}^n$ we write
$$\Gamma(X)=\EE \max_{x \in X} \langle c, x \rangle,$$
where $c=(\gamma_1, \ldots, \gamma_n)$ is a random vector sampled 
from the distribution $\mu^{\otimes n}$ on ${\Bbb R}^n$
and $\langle \cdot, \cdot \rangle$ is the standard scalar product 
in ${\Bbb R}^n$. 
\subhead (4.1) Preliminaries \endsubhead
It is easy to check that 
$$\Gamma(X) \geq \Gamma(Y) \quad \text{provided} \quad 
Y \subset X $$ and that 
$$\Gamma(X)=0 \quad \text{if} \quad |X|=1,$$
that is, if $X$ is a point ($\mu$ is symmetric). 
It follows that $\Gamma(X) \geq 0$ for any finite non-empty subset 
$X \subset {\Bbb R}^n$. Moreover 
$$\Gamma(X+Y)=\Gamma(X)+\Gamma(Y) \quad 
\text{where} \quad X+Y=\bigl\{x+y: \ x \in X,\ y \in Y \bigr\}$$
is the Minkowski sum of $X$ and $Y$. In particular, 
$\Gamma(X+y)=\Gamma(X)$ for any set $X$ and any point $y$.
We note that 
$$\Gamma(\lambda X)=|\lambda| \Gamma(X) \quad 
\text{where} \quad \lambda X=\bigl\{ \lambda x:\ x \in X \bigr\}$$
is a dilation of $X$ and that $\Gamma(X)$ is invariant under the 
action of the hyperoctahedral group, which permutes and changes 
signs of the coordinates. 

Let $S(k, n)$ be the Hamming sphere of radius $k$ centered at 
the origin, that is, the set of points 
$x=(\xi_1, \ldots, \xi_n) \in \{0, 1\}^n$ such that 
$\xi_1 + \ldots + \xi_n=k$. Combinatorially, $S(k, n)$ is the family of 
all $k$-subsets of $\{1, \ldots, n\}$.

Let $F$ be the cumulative distribution function of $\mu$. 
In this section, we prove the following main result.
\proclaim{(4.2) Theorem} For a non-empty set $X \subset \{0, 1\}^n$
and a number $\tau > 0$, let 
$$G(X, \tau)=\ln |X| - \tau \Gamma(X).$$
Let 
$$g_{\tau}(a)=\ln(1+e^{-\tau a}) -  \tau \int_a^{+\infty} 
\bigl(1-F(t)\bigr) \ dt \quad \text{for} \quad 
a \in {\Bbb R}.$$
\roster
\item  For any non-empty set $X \subset \{0, 1\}^n$, we 
have 
$$G(X, \tau) \leq n \sup_{a \geq 0} g_{\tau}(a);$$
\item Suppose that 
$$\sup_{a \geq 0} g_{\tau}(a)=g_{\tau}(a_0) >0 \quad \text{for some,
necessarily finite,} \quad
a_0 \geq 0.$$
Then there exists a sequence $X_n=S(k_n, n) \subset \{0, 1\}^n$ of 
Hamming spheres such that
$$\lim_{n \longrightarrow +\infty} {G(X_n, \tau) \over n} =
g_{\tau}(a_0).$$
Assuming that $F$ is continuous and strictly increasing, we can
choose $k_n=\alpha n +o(n)$ for $\alpha=1-F(a_0)$.
\endroster
\endproclaim
Before we embark on the proof of Theorem 4.2, we summarize 
some useful properties of $g_{\tau}(a)$.
\subhead (4.3) Properties of $g_{\tau}$ \endsubhead
We observe that 
$$g_{\tau}(0)=\ln 2 - \tau \int_0^{+\infty} \bigl(1-F(t)\bigr) \ dt=
\ln 2 -\tau \int_0^{+\infty} t \ dF(t).$$
Furthermore,
$$\lim_{a \longrightarrow +\infty} g_{\tau}(a)=0,$$
since $\mu$ has expectation.
If $F(t)$ is continuous then $g_{\tau}$ is differentiable and 
$$g_{\tau}'(a)=\tau \Bigl( {e^{\tau a} \over 1+e^{\tau a}} 
-F(a) \Bigr).$$
In particular, $a$ is a critical point of $g_{\tau}(a)$ if 
and only if $a$ is a solution of the equation 
$${e^{\tau a} \over 1+ e^{\tau a}}=F(a)$$
or, in other words, if  
$$a\tau=\ln F(a) - \ln \bigl(1-F(a)\bigr).$$
In particular, $a=0$ is always a critical point of $g_{\tau}$.
\bigskip
We prove Part (1) of Theorem 4.2 by induction on $n$, inspired by 
Talagrand's method \cite{T95}. The induction is 
based on the following simple observation.
\proclaim{(4.4) Lemma} 
Suppose that the cumulative distribution function $F$ is 
continuous.
For a non-empty set $X \subset \{0, 1\}^n$,
$n>1$, let 
$$\split &X_1=\Bigl\{x \in \{0,1 \}^{n-1}: \ (x, 1) \in X \Bigr\} 
\quad \text{and} \\ 
&X_0=\Bigl\{x \in \{0, 1\}^{n-1}: \ (x, 0) \in X \Bigr\}.
\endsplit$$  
Then, for any 
$a \in {\Bbb R}$ we have 
$$\Gamma(X) \geq \bigl(1-F(a)\bigr) \Gamma(X_1) + 
F(a) \Gamma(X_0) + \int_a^{+\infty} t \ dF(t).$$
\endproclaim
\demo{Proof} Let $c=(\overline{c}, \gamma)$, where 
$\overline{c} \in {\Bbb R}^{n-1}$,  
$\gamma \in {\Bbb R}$, and let 
$$w(X,c)=\max_{x \in X} \langle x, c \rangle \quad \text{for} \quad 
c \in {\Bbb R}^n.$$
Clearly,
$$w(X, c) \geq w(X_1, \overline{c}) + \gamma \quad 
\text{and} \quad w(X, c) \geq w(X_0, \overline{c}).$$
Therefore, 
$$\split 
\Gamma(X)=&\int_{{\Bbb R}^n} w(X, c) \ d\mu^{\otimes n}(c) \\
= &\int_{{\Bbb R}^n: \gamma > a} w(X, c) \ d \mu^{\otimes n}(c) +
\int_{{\Bbb R}^n: \gamma \leq a} w(X, c) \ d \mu^{\otimes n}(c) \\
\geq &\int_{{\Bbb R}^n: \gamma > a} \Bigl(w(X_1, \overline{c}) + 
\gamma\Bigr) d \mu^{\otimes n}(c) + \int_{{\Bbb R}^n: \gamma \leq a} 
w(X_0, \overline{c}) \  d \mu^{\otimes n}(c) \\
=& \bigl(1-F(a)\bigr) \int_{{\Bbb R}^{n-1}} w(X_1, \overline{c}) 
\ d \mu^{\otimes n-1}(\overline{c}) + \int_a^{+\infty} \gamma \ d F(\gamma)
\\ &\qquad + F(a) \int_{{\Bbb R}^{n-1}} w(X_0, \overline{c}) \ 
d \mu^{\otimes n-1}(\overline{c}) \\ =
& \bigl( 1-F(a)\bigr) \Gamma(X_1) + F(a) \Gamma(X_0) +
 \int_a^{+\infty} \gamma \ d F(\gamma), \endsplit$$
 and the proof follows.
{\hfill \hfill \hfill} \qed
\enddemo
\proclaim{(4.5) Lemma} 
Suppose that the cumulative distribution function $F$ is 
continuous.
For a non-empty set 
$X \subset \{0,1\}^n$ and a number $\tau > 0$ let 
$G(X, \tau)$ and
$g_{\tau}(a)$ be defined as in Theorem 4.2.
Then for any non-empty set $X \subset \{0,1\}^n$, $n>1$,
there exists a non-empty set $Y \subset \{0,1\}^{n-1}$ 
such that 
$$G(X, \tau) \leq G(Y, \tau) + \sup_{a \geq 0} g_{\tau}(a).$$  
\endproclaim
\demo{Proof}
Let us construct $X_1$ and $X_0$ as in Lemma 4.4. We have 
$$|X_1|=\lambda |X| \quad \text{and} \quad |X_0|=(1-\lambda) |X|
\quad  \text{for some} \quad 0 \leq \lambda \leq 1.$$ Without loss of 
generality, we assume 
that $0 \leq \lambda \leq 1/2$. Otherwise, we replace $X$ by $X'$,
where
$$X'=\Bigl\{ (\xi_1, \ldots, 1-\xi_n): \quad (\xi_1, \ldots, \xi_n) \in X
\Bigr\}.$$
Clearly, $|X|=|X'|$ and by Section 4.1, $\Gamma(X)=\Gamma(X')$.

If $\lambda=0$ we choose $Y=X_0$. Identifying ${\Bbb R}^{n-1}$ with 
the hyperplane $\xi_n=0$ in ${\Bbb R}^n$, we observe that $X=Y$ and 
so $G(X,\tau)=G(Y, \tau)$. 
 Since by Section 4.3 
we have 
$$\sup_{a \geq 0} g_{\tau}(a) \geq 0,$$ the result follows.

Thus we assume that $0<\lambda \leq 1/2$. Let 
$Y \in \{ X_0, X_1\}$ be the set with the larger value of
$G(\cdot, \tau)$, where the ties are broken arbitrarily. We have 
$$|X|={1 \over \lambda} |X_1| \quad \text{and} \quad 
|X|={1 \over 1-\lambda} |X_0|.$$
For any $a \geq 0$ 
$$\split G(X, \tau)=&\ln |X| -\tau \Gamma(X)=
\bigl(1-F(a)\bigr) \ln |X| + F(a) \ln |X| 
-\tau \Gamma(X) \\=
&\bigl(1-F(a)\bigr) \ln |X_1| + F(a) \ln |X_0|
\\ &\quad  + \bigl((1-F(a)\bigr) \ln 
{1 \over \lambda} +F(a) \ln {1 \over 1-\lambda} - \tau \Gamma(X).
\endsplit$$
By Lemma 4.4 we conclude that 
$$\split G(X, \tau) \leq & \bigl(1-F(a)\bigr) \ln |X_1| + F(a) \ln |X_0|
\\ &\quad  + \bigl((1-F(a)\bigr) \ln 
{1 \over \lambda} +F(a) \ln {1 \over 1-\lambda} \\
&\quad -\bigl(1-F(a)\bigr) \tau \Gamma(X_1) - F(a) \tau \Gamma(X_0) -
\tau \int_a^{+\infty} t \ dF(t) \\
=&\bigl(1-F(a)\bigr) G(X_1, \tau) +F(a) G(X_0, \tau) \\  
& \quad + \bigl(1-F(a)\bigr) \ln 
{1 \over \lambda} +F(a) \ln {1 \over 1-\lambda} -
\tau \int_a^{+\infty} t \ dF(t)  \\ 
\leq &G(Y, \tau) +\bigl(1-F(a)\bigr) \ln 
{1 \over \lambda} + F(a) \ln {1 \over 1-\lambda} -
\tau \int_a^{+\infty} t\ dF(t).\endsplit$$
Optimizing in $a$, we choose 
$$a={1 \over \tau} \ln \Bigl( {1 -\lambda \over \lambda}\Bigr), \quad
\text{so} \quad a \geq 0.\tag4.5.1$$
Then 
$${1 \over \lambda} = 1+e^{\tau a} \quad \text{and} \quad 
{1 \over 1-\lambda}={1+e^{\tau a} \over e^{\tau a}}.$$
Hence 
$$\split G(X, \tau) \leq &G(Y, \tau)+\ln (1+e^{\tau a}) -\tau a F(a) -
\tau \int_a^{+\infty} t \ dF(t)  \\
=&G(Y, \tau)+\ln(1+e^{-\tau a}) + \tau a \bigl(1-F(a)\bigr) + 
\tau \int_a^{+\infty} t \ d\bigl(1-F(t)\bigr) \\  
= & G(Y, \tau) + \ln (1 + e^{-\tau a}) -
\tau \int_a^{+\infty} \bigl(1-F(t)\bigr)\ dt \\
  = &G(Y, \tau) +g_{\tau}(a),  \endsplit$$
as claimed.
{\hfill \hfill \hfill} \qed 
\enddemo
Now we are ready to prove Part (1) of Theorem 4.2.
\demo{Proof of Part (1) of Theorem 4.2}
Without loss of generality, we may assume 
that the cumulative distribution function $F$ is continuous.
The proof follows by induction on $n$. 
For $n=1$, there are two possibilities. If $|X|=1$ then 
$G(X, \tau)=0$ (see Section 4.1) and the result holds
since 
$$\sup_{a \geq 0} g_{\tau}(a) \geq 0,$$
see Section 4.3.
If $|X|=2$ then $X=\{0,1\}$ and 
$$G(X, \tau)=\ln 2 -\tau \int_0^{+\infty} t \ dF(t)=g_{\tau}(0),$$
so the inequality holds as well.

The induction step follows by Lemma 4.5.
{\hfill \hfill \hfill} \qed
\enddemo

Let $S(k, n)$ be the Hamming sphere of radius $k$, that is, the set of all
$k$-subsets of $\{1, \ldots, n\}$. Given 
weights $\gamma_1, \ldots, \gamma_n$, the maximum weight of a 
subset $x \in S(k,n)$ is the sum of the first $k$ largest weights among 
$\gamma_1, \ldots, \gamma_n$. 

The proof of Part (2) of Theorem 4.2 is based on the following lemma.
\proclaim{(4.6) Lemma} Suppose that the cumulative distribution 
function $F$ of $\mu$ is strictly increasing and continuous. Let us choose 
$0 < \alpha < 1$ and let $X_n$ be the Hamming sphere 
of radius $\alpha n+o(n)$ in $\{0,1\}^n$.

Then 
$$\lim_{n \longrightarrow +\infty} 
{\Gamma(X_n) \over n}=
\int_{F^{-1}(1-\alpha)}^{+\infty} t \ d F(t).$$ 
\endproclaim
\demo{Proof} 
Let $\gamma_1, \ldots, \gamma_n$ be independent 
random variables with the distribution $\mu$ and let 
$u_{1:n} \leq u_{2:n} \leq \ldots \leq u_{n:n}$ be the corresponding 
order statistics, that is, the permutation of 
$\gamma_1, \ldots, \gamma_n$ in the increasing order.
Then
$$\max_{x \in X_n} \sum_{i \in x} \gamma_i= 
\sum_{m=n-\alpha n +o(n)}^n u_{m:n}.$$
Consequently, $\Gamma(X_n)$ is the expectation of the last sum.

The corresponding asymptotics for the order statistics is well known,
see, for example, \cite{S73}.
{\hfill \hfill \hfill} \qed     
\enddemo

Now we are ready to complete the proof of Theorem 4.2.
\demo{Proof of Part (2) of Theorem 4.2}
Without loss of generality, we assume that the cumulative
distribution function $F$ of $\mu$ is continuous and strictly 
increasing. Let us choose $\alpha$ and $k_n$ as 
described, so $X_n \subset \{0,1\}^n$ is 
the Hamming sphere of radius $\alpha n+o(n)$ in $\{0,1\}^n$.

As is known (see, for example, Theorem 1.4.5 of \cite{Li99}),
$$\split \lim_{n \longrightarrow +\infty}
{\ln |X_n| \over n}= &\alpha  \ln {1 \over \alpha} + 
(1-\alpha) \ln {1 \over 1-\alpha}\\=
&\bigl(1-F(a_0)\bigr)\ln {1 \over 1-F(a_0)} 
+F(a_0) \ln {1 \over F(a_0)}. \endsplit$$
Moreover, by Lemma 4.6, 
$$\lim_{n \longrightarrow +\infty} 
{\Gamma(X_n) \over n}= \int_{a_0}^{+\infty} t \ dF(t).$$
Hence
$$\split &\lim_{n \longrightarrow +\infty} 
{G(X_n, \tau) \over n} \\=
&\qquad 
\bigl(1-F(a_0)\bigr) \ln {1 \over 1-F(a_0)} + F(a_0) \ln {1 \over F(a_0)} -
\tau \int_{a_0}^{+\infty} t \ dF(t). \endsplit $$

On the other hand, 
$$\split g_{\tau}(a)=&\ln(1+e^{-\tau a})-\tau 
\int_a^{+\infty} \bigl(1-F(t)\bigr) \ dt \\
= &\ln(1+e^{-\tau a})+\tau a\bigl(1-F(a)\bigr)- \tau \int_a^{+\infty} 
t \ dF(t).\endsplit$$
Since $a_0$ is a critical point of $g_{\tau}$, we have
$$\tau a_0=\ln F(a_0)-\ln \bigl(1-F(a_0)\bigr),$$
cf. Section 4.3.
Therefore,
$$\split g_{\tau}(a_0)=&-\ln F(a_0) +
\Bigl(\ln F(a_0)- \ln \bigl(1-F(a_0)\bigr)\Bigr)
\bigr(1-F(a_0)\bigl) - \tau \int_{a_0}^{+\infty} t \ d F(t)\\
=&\bigl(1-F(a_0)\bigr) \ln {1 \over 1-F(a_0)} +F(a_0) \ln {1 \over F(a_0)}
-\tau \int_{a_0}^{+\infty} t \ dF(t) \endsplit $$
and the proof follows.
{\hfill \hfill \hfill} \qed
\enddemo

Some remarks are in order.
\remark{(4.7) Remarks} 

(4.7.1) Optimizing in $a$ in Lemma 4.4, we substitute 
$a=\Gamma(X_0)-\Gamma(X_1)$
and obtain the inequality
$$\split
\Gamma(X) \geq &\bigl(1-F(a)\bigr)\Gamma(X_1)+F(a)\Gamma(X_0)+ 
\int_a^{+\infty} t \ dF(t) \\=&\Gamma(X_0) +\int_{\Gamma(X_0)-\Gamma(X_1)}
\bigl(1-F(t)\bigr) \ dt. \endsplit$$
This inequality is harder to work with than with 
that of Theorem 4.2 but it sometimes leads to more 
delicate estimates, see Section 7.

(4.7.2) M. Talagrand proved in \cite{T94} that for every non-empty 
set $X$ of subsets of $\{1, \ldots, n\}$ there is a ``shifted'' set 
$X'$ of subsets of $\{1, \ldots, n\}$ such that $|X'|=|X|$, 
$\Gamma(X') \leq \Gamma(X)$, $X'$ is {\it hereditary} 
(that is, if $x \in X'$ and $y \subset x$ then $y \in X'$) and 
{\it left-hereditary} (that is, if $x \in X'$, $i \in x$, $j \notin x$ 
and $j < i$ then the subset $x \cup \{j\} \setminus \{i\}$ also lies
in $X'$). 
\endremark

\head 5. General Estimates: the Lower Bound \endhead

Let us fix a symmetric probability
 measure $\mu$ with the cumulative distribution 
function $F$.
In this section, we prove the following main result.
\proclaim{(5.1) Theorem} Assume that the moment generating function 
$$L(\delta, \mu)=L(\delta)=\EE e^{\delta x}=\int_{-\infty}^{+\infty} 
e^{\delta x} \ d \mu(x)$$
is finite in some neighborhood of $\delta=0$.
Let
$$h(t, \mu)=h(t) =\sup_{\delta \geq 0} \Bigl( \delta t - \ln L(\delta) \Bigr)
\quad \text{for} \quad t \geq 0.$$
\roster
\item For any non-empty family $X$ of $k$-subsets 
of $\{1, \ldots, n\}$, we have 
$$k^{-1} \ln |X| \geq h (t) \quad \text{for} \quad t=k^{-1} \Gamma(X);$$
\item For any $t>0$ such that $F(t)<1$  
and for any $0< \epsilon < 0.1$ there exist
$k=k(t, \epsilon, \mu)$, $n=n(k)$, and 
a family of $k$-subsets of the set 
$\{1, \ldots, n\}$ such that
$$k^{-1} \Gamma(X) \geq (1-\epsilon)t \quad \text{and} \quad 
k^{-1}\ln |X| \leq  h(t) +\epsilon.$$
\endroster
\endproclaim
Before proving Theorem 5.1, we summarize some properties of $L(\delta)$ and $h(t)$.
\subhead (5.2) Preliminaries \endsubhead
Let $f(\delta)=\ln L(\delta)$. Thus we assume that $f(\delta)$ is finite 
on some interval in ${\Bbb R}$, possibly on the whole line.
It is known that $f(\delta)$ is convex and continuous on the 
interval where it is finite, see, for example, Section 5.11 
of \cite{GS01}.
Since $\mu$ is symmetric, we have $f(0)=0$ and 
from Jensen's inequality we conclude that 
$f(\delta) \geq 0$ for all $\delta$.

The function $h(t)$ is convex conjugate to $f(\delta)$.
Therefore, $h(t)$ is finite on some interval where it is 
convex, continuous and approaches $+\infty$ as $t$ approaches a boundary 
point not in the interval. Besides, 
$$h(t)={t^2 \over 2 D} +O(t^4) \quad \text{for} \quad t \approx 0,$$
where $D$ is the variance of $\mu$. In particular, $h(0)=0$ and 
$h(t)$ is increasing for $t \geq 0$, see Section 5.11 of \cite{GS01}.
\bigskip
Now we are ready to prove Theorem 5.1.
\demo{Proof of Theorem 5.1}

Let us prove Part (1).
Without loss of generality, we assume that $\Gamma(X)>0$.
Let us choose a positive integer $m$, let $N=nm$, $K=km$ and let 
$$Y=\underbrace{X \times \ldots \times X}_{\text{$m$ times}} 
\subset \{0, 1\}^N.$$
Let us pick a point $y=(x_1, \ldots, x_m)$ from $Y$, where $x_i \in X$ 
for $i=1, \ldots, m$. Thus some $K$ coordinates of $y$ are 1's and 
the rest are 0's. Let us endow ${\Bbb R}^N$ with 
the product measure 
$\mu^{\otimes N}$ and let $\gamma_1, \ldots, \gamma_K$ be independent 
random variables with the distribution $\mu$.
Then, for any $t>0$
$$\PP\Bigl\{ c \in {\Bbb R}^N: \ \langle c, y \rangle > m t \Bigr\}=
\PP \Bigl\{ \sum_{i=1}^K \gamma_i > m t\Bigr\}=
\PP \Bigl\{ \sum_{i=1}^K \gamma_i >  K{t \over k} \Bigr\}.$$ 
By the Large Deviations Inequality (see, for example, Section 5.11 of 
\cite{GS01}) 
$$\PP \Bigl\{ \sum_{i=1}^K \gamma_i \geq K{t \over k} \Bigr\} \leq 
\exp\bigl\{-K h(t/k)\bigr\}.$$
Therefore,
$$\PP \Bigl\{ c \in {\Bbb R}^N: \quad \max_{y \in Y} \langle c, y \rangle 
> mt \Bigr\} \leq |Y| \exp\bigl\{-K h(t/k) \bigr\}=
\Bigl(|X| \exp\bigl\{-k h(t/k) \bigr\} \Bigr)^m.$$
Since a vector $c \in {\Bbb R}^N$ is an $m$-tuple $c=(c_1, \ldots, c_m)$
with $c_i \in {\Bbb R}^n$ and
$$\max_{y \in Y} \langle c, y \rangle = \sum_{i=1}^m \max_{x \in X} 
\langle c_i, x \rangle,$$ 
the last inequality can be written as 
$$\PP\Bigl\{ c_1, \ldots, c_m: \quad 
{1 \over m} \sum_{i=1}^m \max_{x \in X} \langle c_i, x \rangle > t \Bigr\} 
\leq \Bigl(|X| \exp\bigl\{-k h(t/k)\bigr\} \Bigr)^m.$$
However, by the Law of Large Numbers
$${1 \over m} \sum_{i=1}^m \max_{x \in X} \langle c_i, x \rangle
\longrightarrow \Gamma(X) \quad \text{in probability}$$ 
as $m \longrightarrow +\infty$. Therefore, for any $0<t< \Gamma(X)$, 
$$\PP\Bigl\{ c_1, \ldots, c_m: \quad 
{1 \over m} \sum_{i=1}^m \max_{x \in X} \langle c_i, x \rangle > t \Bigr\} 
\longrightarrow 1 \quad \text{as} \quad m \longrightarrow +\infty.$$
Therefore, we must have
$$|X| \exp\{-k h(t/k) \}  \geq 1 \quad \text{for every} \quad t< \Gamma(X).$$
Hence
$$k^{-1} \ln |X| \geq h(t) \quad \text{for every} \quad t < k^{-1} \Gamma(X),$$
and the proof follows by the continuity of $h$, cf. Section 5.2.

Let us prove Part (2).
Let $\gamma_1, \ldots, \gamma_k$ be independent random 
variables having the distribution $\mu$. By the Large Deviations 
Theorem (see Section 5.11 of \cite{GS01}), if $k=k(\epsilon,t,\mu)$ is 
sufficiently large then 
$$\PP \Bigl\{ \sum_{i=1}^k \gamma_i > k t \Bigr\} \geq 
\exp\bigl\{-k \bigl(h(t) +\epsilon/2 \bigr)\bigr\}.$$  
We make $k$ large enough to ensure, additionally, that 
$\bigl(\ln 3 + \ln \ln (1 / \epsilon) \bigr)/k \leq \epsilon/2$.

Let $|X|$ be the largest integer not exceeding 
$$ 3\ln {1 \over \epsilon}\exp\bigl\{ k \bigl(h(t)+\epsilon/2 \bigr) \bigr\},$$
so $k^{-1} \ln |X| \leq h(t)+\epsilon$,
and let $X$ consist of $|X|$ pairwise disjoint $k$-subsets of 
$\{1, \ldots, n\}$ for a sufficiently large $n=n(k)$.

Suppose that $c=(\gamma_1, \ldots, \gamma_n)$ is a random vector 
of independent weights with the distribution $\mu$. 
Since $x \in X$ are disjoint,
the weights $\sum_{i \in x} \gamma_i$ 
of subsets from $X$ are independent random variables. Let $w(X,c)$ 
be the largest weight of a subset $x \in X$. We have
$$\PP\Bigl\{c: \quad w(X, c) \leq  kt \Bigr\} \leq 
\Bigl(1-\exp\bigl\{-k \bigl(h(t) +\epsilon/2 \bigr)\bigr\} \Bigr)^{|X|}
\leq \epsilon/2.$$
Similarly (since $\mu$ is symmetric):
$$\PP\Bigl\{c: \quad w(X, -c) \leq  kt \Bigr\} \leq \epsilon/2,$$
and, therefore,
$$\PP\Bigl\{c: \quad w(X, c) +w(X,-c)  \leq  2kt \Bigr\} \leq \epsilon.$$
Since $w(X, c) + w(X, -c)$ is always non-negative, its 
expectation is at least $(1-\epsilon)2kt$. On the other hand, 
this expectation is $2 \Gamma(X)$. Hence we have constructed 
a family $X$ of $k$-subsets such that
$$k^{-1}\Gamma(X) \geq (1-\epsilon)t \quad \text{and} \quad 
k^{-1} \ln |X| \leq h(t) +\epsilon.$$ 
{\hfill \hfill \hfill} \qed 
\enddemo
\remark{(5.3) Remarks} 

(5.3.1) Using the convexity of $h(t)$, one can extend 
the bound of Part (1) of Theorem 5.1 to families $X$ of {\it at most} $k$-element
subsets of $\{1, \ldots, n\}$.

(5.3.2) Suppose that the moment generating function 
$L(\delta, \mu)$ is infinite for all $\delta$ except for $\delta=0$. 
Let us choose $t>0$ and $0 < \epsilon < 0.1$.
We claim that there exists 
a family $X$ of $k$-subsets of $\{1, \ldots, n\}$ such that
$$k^{-1} \Gamma(X) \geq (1-\epsilon)t \quad \text{and} \quad 
k^{-1} \ln |X| \leq \epsilon$$
(in other words, we can formally take $h(t) \equiv 0$ in Part (2) 
of Theorem 5.1). Let $\gamma$ be a random variable with the distribution $\mu$. For $c>0$, let $\gamma_c$ be the {\it truncation} of $\gamma$:
$$\gamma_c=\cases \gamma, &\text{if} \quad |\gamma| \leq c \\
0,  &\text{if} \quad |\gamma|>c. \endcases$$
Let $\mu_c$ be the distribution of $\gamma_c$. It is not hard to see 
that $\Gamma(X, \mu) \geq \Gamma(X, \mu_c)$
(consider $\Gamma(X)$ as the expectation of
$0.5w(X,c)+0.5w(X,-c)$, where $w(X,c)$ is the maximum weight of a subset 
$x \in X$ for the vector $c=(\gamma_1, \ldots, \gamma_n)$ of weights).
Choosing a sufficiently 
large $c$ brings $h(t, \mu_c)$ arbitrarily close to 0.
Then we construct a set $X$ as in Part (2) of Theorem 5.1. 

(5.3.3) Our proof of Part (2) of Theorem 5.1 seems to require $n$ to
be exponentially large in $k$. This is not so, since every suitable 
pair $n,k$ can be rescaled to a suitable pair $N=nm$, $K=km$ for 
a positive integer $m$.
Let $X$ be a family of $k$-subsets of $\{1, \ldots, n\}$
constructed in the proof of Part (2) and
let  $$Y=\underbrace{X \times \ldots \times X}_{\text{$m$ times}} 
\subset \{0, 1\}^N.$$
Then $Y$ is a family of $K$-subsets of $\{1, \ldots, N\}$ and
$$K^{-1} \Gamma(Y) \geq (1-\epsilon)t \quad \text{and} \quad 
K^{-1}\ln |Y| \leq  h(t) +\epsilon.$$
\endremark

\head 6. The Logistic Measure: Proofs \endhead

In this section, we prove Theorems 3.1 and 3.3. 
\demo{Proof of Theorem 3.1} To prove Part (1), let us choose
$\tau=1$ in Part (1) of Theorem 4.2.
We have 
$$g_1(a)=\ln(1+e^{-a})-\int_a^{+\infty} {e^{-t} \over 1+e^{-t}}\ dt=
0 \quad 
\text{for all} \quad a.$$
Hence 
$$\ln |X| \leq \Gamma(X)$$
as claimed.

To prove Part (2), we use Part (1) of Theorem 5.1. The moment generating 
function of the logistic distribution is given by 
$$L(\delta)=\int_{-\infty}^{+\infty} 
{e^{\delta x} \over e^x +e^{-x} +2} \ dx=
{\pi \delta \over \sin \pi \delta} \quad \text{for} \quad -1 < \delta <1,$$
see \cite{M85}.
Hence the formula for $h(t)$ follows.

It follows from Section 5.2 that $h$ is convex and increasing.
{\hfill \hfill \hfill} \qed
\enddemo

Now we are ready to prove optimality of the logistic distribution. 

\demo{Proof of Theorem 3.3} 
Let us choose $\mu \in {\Cal M}$ and let $F_{\mu}$ be the 
cumulative distribution function of $\mu$. We claim that $F_{\mu}(t)<1$ 
for all $t \in {\Bbb R}$. To see that, we let $\tau=1$ in Theorem 4.2.
If $F_{\mu}(t)=1$ then $g_1(t)>0$ and, by Part (2) 
of Theorem 4.2, there is a set $X \subset \{0, 1\}^n$ with $\ln |X| > \Gamma(X)$, which contradicts the definition of ${\Cal M}$.

 Let us assume first that the 
moment generating function $L(\delta, \mu)$ is finite in some 
neighborhood of $\delta=0$. Then, by Theorem 5.1, 
we have $c(t, \mu)=h(t, \mu)$ and hence we must prove that 
$h(t, \mu) \leq h(t, \mu_0)$, where $\mu_0$ is the logistic distribution.

Let 
$$T(a, \mu)=\int_a^{+\infty} \bigl(1-F_{\mu}(t)\bigr) \ dt.$$
We can write
$$\split &\int_0^{+\infty} e^{\delta x} \ d F_{\mu}(x)=
-\int_0^{+\infty} e^{\delta x} \ d \bigl(1-F_{\mu}(x)\bigr)=
{1 \over 2} + \int_0^{+\infty} \delta e^{\delta x} 
\bigl(1-F_{\mu}(x)\bigr) 
\ dx \\ = &{1 \over 2} + \int_0^{+\infty} \delta e^{\delta x}  \ 
d(-T(x, \mu))=
{1 \over 2} + \delta T(0, \mu) + \int_0^{+\infty} \delta^2 e^{\delta x} 
T(x, \mu) \ dx. \endsplit$$
Similarly,
$$\int_{-\infty}^{0} e^{\delta x} \ d F_{\mu}(x)=
\int_0^{+\infty} e^{-\delta x} \ d F_{\mu}(x) =
{1 \over 2} - \delta T(0, \mu) + \int_0^{+\infty} \delta^2 e^{-\delta x} 
T(x, \mu) \ dx.$$
Therefore,
$$L(\delta, \mu)=1 + \delta^2 \int_0^{+\infty} 
\bigl(e^{-\delta x} + e^{\delta x} \bigr) T(x, \mu) \ dx.$$ 

Since $\ln |X| \leq \Gamma(X)$, by Part (2) of Theorem 4.2 we conclude 
that
$$T(a, \mu) \geq \ln(1+e^{-a})=T(a, \mu_0) \quad \text{for all} 
\quad a \geq 0.$$
Therefore, $L(\delta, \mu) \geq L(\delta, \mu_0)$ and
$h(t, \mu) \leq h(t, \mu_0)$ for all $t \geq 0$, as claimed.

Suppose now that the moment generating function $L(\delta, \mu)$ is 
infinite for $\delta \ne 0$. Then, as follows from Remark 5.3.2, 
$c(t, \mu)=0$ for all $t>0$, which completes the proof.
{\hfill \hfill \hfill} \qed
\enddemo

\head 7. The Exponential Measure \endhead

Let us choose $\mu$ to be the measure with
density 
$${1 \over 2} e^{-|\gamma|} \quad \text{for} \quad 
\gamma \in {\Bbb R}.$$
As we have already mentioned, one of the results of \cite{T94} 
is the estimate 
$$\ln |X| \leq c \Gamma(X)$$ for some absolute 
constant $c$. In this section, we find the optimal value of $c$ and 
establish some general isoperimetric inequalities which, we believe,
are interesting in their own right.
\proclaim{(7.1) Theorem} Let $\mu$ be the measure with density 
$e^{-|\gamma|}/2$ for $\gamma \in {\Bbb R}$.
\roster 
\item
Let $X \subset \{0,1\}^n$ be a non-empty subset of the Boolean cube.
Then
$$\ln |X| \leq (2 \ln 2) \Gamma(X);$$
\item Let $X \subset \{0,1\}^n$ be a non-empty subset 
of the Boolean cube such that $\xi_1 + \ldots + \xi_n \leq k$ 
for every $(\xi_1, \ldots, \xi_n) \in X$. That is, $X$ lies in the 
Hamming ball of radius $k$ and we may interpret $X$ as a family of 
at most $k$-element subsets of $\{1, \ldots, n\}$.
Then 
$$\ln |X| \leq \Gamma(X) + k \ln 2.$$
\endroster
\endproclaim
Before we prove Theorem 7.1, we note that 
$c=2\ln 2$ is the best possible 
value in Part (1).
If $X$ is a $m$-dimensional face of the Boolean cube 
then 
$\ln |X|=m \ln 2$ and we show that $\Gamma(X)=m/2$, so the 
equality holds. As in Section 3, it suffices to check the formula for 
$X=\{0,1\}$, in
which case
$$\Gamma(X)={1 \over 2} \int_0^{+\infty} x e^{-x} \ dx={1 \over 2}.$$ 
The inequality of Part (2) is asymptotically sharp: if $X$ is the Hamming 
sphere of radius $k=o(n)$ in $\{0,1\}^n$, then 
$$\Gamma(X)=\ln |X|-k \ln 2 +o(k) \quad \text{as} \quad 
k \longrightarrow +\infty,$$
cf. Lemma 4.6. 

As for the lower bound, using Part (1) of Theorem 5.1 
one can show that for any 
non-empty family $X$ of $k$-subsets of $\{1, \ldots, n\}$, we have
$$k^{-1} \ln |X| \geq h\bigl(k^{-1} \Gamma(X)\bigr),$$ 
where 
$$\split h(t)=&\sqrt{1 +t^2} + 
\ln \bigl(\sqrt{1 +t^2} -1 \bigr) - 2 \ln t + \ln 2 -1 \\
=&t-\ln t -O(1)  \quad \text{for large} \quad t. \endsplit$$
Thus the exponential distribution also allows us to estimate 
$\ln |X|$ up to a constant factor. However, the estimates are  
not as good as for the logistic distribution.

\demo{Proof of Theorem 7.1} To prove Part (1), we use 
Part (1) of Theorem 4.2.

The function $g_{\tau}(a)$ is given by 
$$g_{\tau}(a)=\ln(1+e^{-\tau a}) -{\tau \over 2} e^{-a} \quad 
\text{for} \quad a \in {\Bbb R}.$$
Let us consider the critical points of $g_{\tau}$. 

We have
$$g_{\tau}'(a)={\tau \over 2}\Bigl({e^{(\tau -1)a} + e^{-a} -2 \over 
1+e^{\tau a}}\Bigr).$$
Since the numerator of the fraction is a linear combination of 
two exponential functions and a constant, it can have at most two real zeros.
We observe that $a=0$ is a zero and that $g_{\tau}'(a)<0$ for small
$a>0$ provided $\tau<2$.

Hence for $\tau < 2$ the function $g_{\tau}$ has at most one 
critical point $a>0$ which has to be a point of local minimum.

Therefore
$$\sup_{a \geq 0} g_{\tau}(a)=\max \bigl\{g_{\tau}(0),\ 0 \bigr\} 
\quad \text{for all} \quad \tau<2.$$
Let us choose $\tau =2\ln 2$. 
Then $g_{\tau}(0)=0$ and we conclude that
$$\sup_{a \geq 0} g_{\tau}(a)=0.$$
By Part (1) of Theorem 4.2, we conclude that 
$$\ln |X| \leq \tau \Gamma(X)=(2 \ln 2) \Gamma(X).$$ 

We prove Part (2) by induction on $n$. If $n=1$, there are two 
cases. If $X$ consists of a single point then $\Gamma(X)=0$, $\ln |X|=0$ 
and the inequality is satisfied. If $X=\{0,1\}$ then $k=1$ and 
$\Gamma(X)=1/2$, hence the inequality holds as well.

Suppose that $n>1$. Clearly, we can assume that $k>0$.
 Without loss of generality, we may assume that   
$X$ is hereditary, see Remark 4.7.2. Let us construct sets 
$X_0, X_1 \subset \{0,1 \}^{n-1}$ as in Lemma 4.4. We note that 
$X_0$ lies in the Hamming ball of radius $k$ and $X_1$ lies in the 
Hamming ball of radius $k-1$. Since $X$ is hereditary, $X_1 \subset X_0$.
Therefore,
$$|X_0| \geq |X_1| \quad \text{and} \quad \Gamma(X_0) \geq \Gamma(X_1).$$ 
The inequality of Remark 4.7.1 gives us
$$\Gamma(X) \geq \Gamma(X_0) + {1 \over 2} 
\exp\bigl\{\Gamma(X_1)-\Gamma(X_0) \bigr\}.$$
Let us consider a function 
$$f(a,b)=a+{1 \over 2}e^{b-a}.$$
It is easy to see that for every $a$ the function is increasing 
in $b$ and that for every $b$ it is increasing on the interval 
$a \geq b-\ln 2$. 

Applying the induction hypothesis to $X_0$ and $X_1$, we conclude
$$\split f\bigl(\Gamma(X_0),\ \Gamma(X_1) \bigr) \geq  
& f\bigl(\Gamma(X_0),\ \ln |X_1|-(k-1) \ln 2 \bigr) \\ \geq 
&f\bigl(\ln |X_0|-k \ln 2,\ \ln |X_1| - (k-1) \ln 2 \bigr).\endsplit$$ 

Therefore,
$$\split \Gamma(X) \geq &\ln |X_0|- k\ln 2 + {|X_1| \over |X_0|}=
\Bigl(\ln |X|-k \ln 2\Bigr) + 
\Bigl({|X_1| \over |X_0|}-\ln {|X_1|+|X_0| \over |X_0|}\Bigr)
\\ = &\Bigl(\ln |X|-k \ln 2\Bigr) + 
\bigl(t-\ln (1+ t)\bigr) \quad \text{for} \quad t={|X_1| \over |X_0|} \\
 \geq &\ln |X| -k \ln 2.  \endsplit$$
The proof now follows.
{\hfill \hfill \hfill} \qed
\enddemo

\head 8. An Asymptotic Solution to the Isoperimetric Problem \endhead

In this section, we discuss what sets 
$X_n \subset \{0,1\}^n$ with 
the smallest ratio $\Gamma(X_n, \mu)/\ln |X_n|$ may look like.
We claim that for any symmetric probability measure
$\mu$ with finite variance and for a sufficiently large 
$n$ we can choose $X_n$ to be the product of at most two Hamming spheres. 
\proclaim{(8.1) Theorem} Let us fix a symmetric probability measure 
$\mu$ and a number 
$$0< \alpha < \ln 2.$$ Then there exist numbers $\beta_i, \lambda_i$,
$i=1,2$, depending on $\alpha$ and $\mu$ only, such that
$$\split &0 \leq \beta_i \leq \lambda_i \quad \text{for} \quad i=1,2, \\
&\lambda_1 +\lambda_2=1 \endsplit$$
and the following holds.

Let $S^i_n$ be 
the Hamming sphere of radius $\beta_i n +o(n)$ in the Boolean cube 
of dimension $\lambda_i n+o(n)$, $i=1, 2$, and let 
$Y_n=S^1_n \times S^2_n$ be the direct product of the spheres considered
as a subset of the Boolean cube of dimension $n$. 

Then 
$$\ln |Y_n|=\alpha n +o(n)$$
and for any sequence of sets $X_n \subset \{0,1\}^n$ such that 
$$\ln |X_n|=\alpha n +o(n),$$
we have 
$$\Gamma(Y_n, \mu) \leq \Gamma(X_n, \mu) +o(n).$$
\endproclaim
\demo{Proof} Let $F$ be the cumulative distribution function of $\mu$.
Without loss of generality, we assume that $F$ is continuous and 
strictly increasing. 
Given $\mu$ and $\alpha$, let us consider the function
$$H(\tau, x)={\alpha \over \tau} - {\ln\bigl(1+e^{-\tau x}\bigr) \over \tau}+
\int_x^{+\infty} \bigl(1-F(t)\bigr) \ dt$$
of two variables $\tau >0$ and $x \geq 0$.

By Part (1) of Theorem 4.2, for any $\tau >0$,
$$n^{-1} \Gamma(X_n) \geq  \inf_{x \geq 0} H(\tau, x) \quad 
\text{provided} \quad \ln |X_n|=\alpha n + o(n). \tag8.1.1$$
We claim that there exists $0 < \tau_0 < +\infty$ such that 
$$\inf_{x \geq 0} H(\tau_0, x) \geq \inf_{x \geq 0} H(\tau, x) \quad
\text{for all} \quad \tau \geq 0.$$ Indeed, since $\alpha < \ln 2$,
$$\inf_{x \geq 0} H(\tau, x) \longrightarrow -\infty
\quad \text{as} \quad \tau \longrightarrow 0+.$$
Also,
$$\inf_{x \geq 0} H(\tau, x) \longrightarrow 0 \quad \text{as} \quad
\tau \longrightarrow +\infty.$$
On the other hand, choosing $x_1>0$ such that 
$$\int_{x_1}^{+\infty} \bigl(1-F(t)\bigr) \ dt =\delta >0$$
and $\tau_1$ such that 
$$\alpha -\ln\bigl(1+e^{-\tau_1 x_1}\bigr) >0 
\quad \text{and} \quad \tau_1^{-1}|\alpha - \ln 2|< \delta $$ we observe that
$$\inf_{x \geq 0} H(\tau_1, x)  >0,$$
which implies that there exists $0 < \tau_0 < +\infty$ maximizing
$\inf_{x \geq 0} H(\tau,x)$.

Our next goal is to show that one can find 
$0 \leq x_1, x_2 \leq +\infty$ such that  
$$H(\tau_0, x_1)=H(\tau_0, x_2)=\inf_{x \geq 0} H(\tau_0,x) \tag8.1.2$$
and such that 
$$\aligned 
&{\tau_0 x_1 \over e^{\tau_0 x_1} +1} + \ln\bigl(1+e^{-\tau_0 x_1}\bigr)
 \geq \alpha \quad \text{and} \\
&{\tau_0 x_2 \over e^{\tau_0 x_2} +1} + \ln\bigl(1+e^{-\tau_0 x_2}\bigr)
\leq \alpha. \endaligned \tag8.1.3$$
(it is possible that $x_1=x_2$ or that $x_2=+\infty$).

For $\epsilon$ in a small neighborhood of $0$, we define 
$x_{\epsilon} \geq 0$ as a point such that 
$$H\bigl(\tau_0+\epsilon,\ x_{\epsilon}\bigr)=\inf_{x \geq 0} 
H(\tau_0+\epsilon, x)$$
(possibly $x_{\epsilon}=+\infty$).
We obtain $x_1$ as a limit point of $x_{\epsilon}$ as 
$\epsilon \longrightarrow 0-$ and $x_2$ as a limit point of
$x_{\epsilon}$ as $\epsilon \longrightarrow 0+$. Clearly, (8.1.2) holds
and it remains to show that (8.1.3) holds as well.

Indeed,
$$\split 
H(\tau_0, x_i) \geq &H\bigl(\tau_0+\epsilon, \ x_{\epsilon}\bigr)=
H(\tau_0, \ x_{\epsilon})+ \epsilon {\partial \over \partial \tau} 
H(\overline{\tau}, \ x_{\epsilon} \bigr) \\ \geq 
& H(\tau_0, x_i)+ \epsilon {\partial \over \partial \tau} 
H(\overline{\tau}, \ x_{\epsilon} \bigr) \endsplit$$
for some $\overline{\tau}$ between $\tau_0$ and $\tau_0+\epsilon$ and 
$i=1,2$.

Besides,   
$${\partial \over \partial \tau} H(\tau, x)=
{1 \over \tau^2} 
\Bigl({\tau x \over 1+e^{\tau x}}+
\ln\bigl(1+e^{-\tau x}\bigr) -\alpha \Bigr),
$$
from which we deduce (8.1.3).

Additionally, from (8.1.2) we deduce that if $0<x_i <+\infty$, we must have
$${\partial \over \partial x} H(\tau_0,x_i)=0,$$
that is, 
$${1 \over e^{\tau_0 x_i} +1} =1-F(x_i), \quad \text{for} \quad 
i=1,2, \tag8.1.4$$
which also holds for $x_i=0$ and $x_i=+\infty$.

Now we are ready to define $\lambda_i$ and $\beta_i$. Namely,
we write 
$$\alpha= \sum_{i=1,2}\lambda_i \Bigl({\tau_0 x_i \over e^{\tau_0 x_i} +1}+
\ln\bigl(1+e^{-\tau_0 x_i}\bigr) \Bigr) 
\quad \text{where} \quad 
\lambda_1, \lambda_2 \geq 0 \quad \text{and} \quad 
\lambda_1 + \lambda_2=1,$$
cf. (8.1.3).

Next, we define $\beta_1$ and $\beta_2$ by 
$$\beta_i={\lambda_i \over e^{\tau_0 x_i} +1} \quad \text{for} \quad 
i=1,2.$$

Let $S^i_n$ be the Hamming sphere of dimension $\lambda_i n +o(n)$ and 
radius $\beta_i n +o(n)$. Using Theorem 1.4.5 of \cite{Li99}, we 
obtain
$$\split {1 \over \lambda_i n}\ln |S^i_n| = &{1 \over e^{\tau_0 x_i} +1} 
\ln \bigl(e^{\tau_0 x_i} +1 \bigr) + 
{e^{\tau_0 x_i} \over e^{\tau_0 x_i} +1} 
\ln \bigl(1+ e^{-\tau_0 x_i} \bigr) +o(1) \\
= &{\tau_0 x_i \over e^{\tau_0 x_i} +1} +\ln\bigl(1+e^{-\tau_0 x_i} \bigr)+
o(1).
\endsplit$$
Thus for $Y_n=S^1_n \times S^2_n$, we have
$$\ln |Y_n|=\alpha n +o(n),$$
as claimed.

By Part (2) of Theorem 4.2 and (8.1.4)
$${\Gamma(S^i_n) \over \lambda_i n}=H(\tau_0, x_i) +o(1).$$
Using (8.1.2) we conclude that for $Y_n=S_n^1 \times S_n^2$, we have 
$$\split {\Gamma(Y_n) \over n}=&{\Gamma(S_n^1) +\Gamma(S_n^2) \over n}=
\lambda_1 H(\tau_0, x_1) + \lambda_2 H(\tau_0, x_2) +o(1) \\=
&\inf_{x \geq 0} H(\tau_0, x)+o(1)  \endsplit.$$

Hence, by (8.1.1), 
$$\Gamma(Y_n) \leq \Gamma(X_n) +o(n),$$
which completes the proof.
{\hfill \hfill \hfill} \qed
\enddemo

\head Acknowledgment \endhead

The authors are very grateful to Nathan Linial for many helpful
discussions and encouragement, in particular, for his suggestion to look
for an optimal solution to the ``inverse isoperimetric problem''.

\head References \endhead 
\refstyle{A}

\enddocument
\end